\documentclass[11pt]{article}
\usepackage[latin1]{inputenc}
\usepackage{amsmath,amsthm,amssymb}
\newtheorem{teo}{Theorem}

\newtheorem{prop}{Proposition}

\newtheorem{lemma}{Lemma}

\def\proof{{\it Proof.}\ }

\def\eq#1{(\ref{#1})}

\def\neweq#1{\begin{equation}\label{#1}}
\def\endeq{\end{equation}}
\def\weak{\rightharpoonup}
\def\ep{\varepsilon}
\def\la{\lambda}

\def\phi{\varphi}
\def\RR{{\mathbb R} }

\def\supp{{\rm supp}}
\def\di{\displaystyle}
\def\ri{\rightarrow}

\title{\sc Singular elliptic problems with
lack of compactness}
\author{Marius GHERGU and
Vicen\c tiu R\u ADULESCU\thanks{Correspondence address: Vicen\c
tiu R\u adulescu, Department of Mathematics, University of
Craiova, 200585 Craiova, Romania, fax: +40-251.41.16.88. E-mail:
{\tt radulescu@inf.ucv.ro}}\\
 \small Departament of Mathematics,
        University of Craiova,
                200585 Craiova, Romania}
\date{}
\begin{document}
\baselineskip12pt

\maketitle

\noindent{\small {\bf Abstract}. We consider the following
nonlinear singular elliptic equation
$$-\mbox{div}\,(|x|^{-2a}\nabla u)=K(x)|x|^{-bp}|u|^{p-2}u+\la g(x)
\quad\mbox{in}\,\,\RR^N,$$ where
$g$ belongs to an appropriate weighted Sobolev space,
 and $p$ 
denotes the  Caffarelli--Kohn--Nirenberg
critical exponent associated to $a$, $b$, and $N$. 
Under some natural assumptions on the positive potential $K(x)$ 
we establish the existence of some
$\la_0>0$ such that the above problem has at least two
distinct solutions provided that $\la\in(0,\la_0)$. The proof
relies on Ekeland's Variational Principle and on the
Mountain Pass Theorem without the Palais-Smale condition,
combined with a weighted
variant of the Brezis-Lieb Lemma.



\section{Introduction and the main result}
Many papers have been devoted in the last decades 
to the study of several questions
concerning degenerate elliptic problems. We start with the
following example
\neweq{unu}
 \left\{\begin{tabular}{ll}
$\mbox{div}\,(a(x)\nabla u)+f(u)=0$ \quad & ${\rm in}\,\ \Omega\,$\\
$u=0$ \quad & ${\rm on}\,\ \partial\Omega\,,$\\
\end{tabular} \right.
\endeq
\noindent where $\Omega$ is an arbitrary domain in $\RR ^N$
($N \geq 1$), and $a$ is a nonnegative function
that may have ``essential" zeroes at some points or even may be unbounded. 
The continuous function $g$ satisfies
$f(0)=0$ and $tf(t)$ behaves like $|t|^p$  as $|t|\rightarrow \infty$, with
$2<p<2^*$, where $2^*$ denotes the critical Sobolev exponent.
Notice that equations of this type come from the consideration of standing 
waves in anisotropic Schr\"odinger equations (see \cite{badiale, sirakov, strauss, wangx}). 
Equations like \eq{unu}
are also introduced as models for several physical phenomena related to 
equilibrium of anisotropic media which possibly are somewhere 
{\it ``perfect" insulators} or {\it ``perfect" conductors} 
(see \cite{dl}, p.~79). Problem \eq{unu} has also some
interest in the framework of optimization and $G$--convergence
(see, e.g., \cite{franchi} and the references therein).

Classical results (see \cite{ar, r}) ensure the existence
and the multiplicity of positive or nodal solutions for problem
\eq{unu}, provided that the differential operator 
$Tu:=\mbox{div}\,(a(x)\nabla u)$ is uniformly elliptic.
Several difficulties occur both in the degenerate case
(if $\inf\limits_{\Omega}a=0$) and in the 
singular case (if $\sup\limits_{\Omega}\,a=+\infty$). 
In these situations the classical methods fail to
be applied directly so that the existence and the multiplicity
results (which hold in the nondegenerate case) may become a
delicate matter that is closely related to some phenomena due to
the degenerate character of the differential equation. These problems
have been intensively studied starting with the pioneering paper by Murthy and
Stampacchia \cite{ms} (we also refer to \cite
{cirmi, fks, p}, as well as to the monograph \cite{stre}).

A natural question that arises in concrete applications
is to see what happens if these elliptic (degenerate
or nondegenerate) problems are affected by a certain perturbation.
It is worth pointing out here that the idea of using perturbation
methods in the treatment of
nonlinear boundary value problems
was introduced by Struwe \cite{stru}. Recently, many authors have
been interested in this kind of perturbation problems involving
both critical and sub- or super-critical Sobolev exponent (see, e.g.,
\cite{cr, rs, t}).
\medskip

Our aim in this paper is to study the following degenerate
perturbed problem
\neweq{doi}
-\mbox{div}\,(|x|^{-2a}\nabla u)=K(x)|x|^{-bp}|u|^{p-2}u+\la
g(x)\, \quad\mbox{in}\,\,\RR^N,
\endeq
where
\neweq{trei}\left\{\begin{array}{ll}
&\di \mbox{for}\ \,N\geq 3:\quad-\infty<a<\frac{N-2}{2}, \quad a< b< a+1,
\quad\mbox{and}\quad p=\frac{2N}{N-2+2(b-a)}\,;\\
&\di \mbox{for}\ \,N=2:\quad-\infty<a<0, \quad a< b< a+1,
\quad\mbox{and}\quad p=\frac{2}{b-a}\,;\\
&\di \mbox{for}\ \,N=1:\quad-\infty<a<-\frac{1}{2}, \quad a+\frac 12< b< a+1,
\quad\mbox{and}\quad p=\frac{2}{-1+2(b-a)}\,.
\end{array}\right.
\endeq

Equation \eq{doi} contains the critical Caffarelli-Kohn-Nirenberg exponent $p$,
defined as in \eq{trei}.
In this critical case, some concentration phenomena may occur, due to the action of 
the noncompact group of dilations in $\RR^N$. The lack of compactness of problem
\eq{doi} is also given by the fact that we are looking for entire solutions, that is,
solutions defined on the whole space.

The reason for which we choose the parameters $a$, $b$, and $p$ to 
satisfy the assumption \eq{trei} has to do with the 
following inequality, due to
 Caffarelli, Kohn, and Nirenberg (see
\cite{ckn}):
\neweq{patru}
\di\left(\,\int_{\RR^N}|x|^{-bp}\,|u|^pdx\right)^{1\slash p}\leq
C_{a,b} \left(\,\int_{\RR^N}|x|^{-2a}|\nabla
u|^2dx\right)^{1\slash 2}\,,
\endeq
for all $u\in C^{\infty}_0(\RR^N)$, where $a$, $b$ and $p$
satisfy the condition \eq{trei}. 
We point out that the inequality \eq{patru} also holds true for $b=a+1$ 
(if $N\geq 1$) and $b=a$ (if $N\geq 3$) but, in these cases, 
the best Sobolev constant $C_{a,b}$ in \eq{patru} is never achieved (see \cite{cw} for details).    
The Caffarelli--Kohn--Nirenberg inequality \eq{patru} contains as particular cases 
the classical
Sobolev inequality (if $a=b=0$) and the Hardy inequality (if $a=0$ and $b=1$); we refer
to \cite{bremar, davies, rsw} for further details.

The extremal functions for \eq{patru} are ground state solutions of the singular
Euler equation
$$-\mbox{div}\, (|x|^{-2a}\,\nabla u)=|x|^{-bp}\, |u|^{p-2}u,\qquad
\mbox{in}\ \, \RR^N\,.$$ This equation has been recently studied (see \cite{cw, ww})
in connection with a complete understanding of the best constants,
the qualitative properties of extremal functions, the existence (or nonexistence)
of minimizers, and the symmetry properties of minimizers.
 
 The function $K$ is
assumed to fulfill

\medskip
\noindent $(K1) \qquad K\in L^{\infty}(\RR^N)$,

\medskip
\noindent $(K2)$\qquad esslim$ _{|x|\ri 0}K(x) =$
esslim$_{|x|\ri\infty}K(x)=K_0\in (0,\infty )$ and $K(x)\geq K_0$
a.e. in $\RR^N$,

\medskip
\noindent$(K3)$\qquad  meas$\di\left(\,\{x\in
\RR^N:K(x)>K_0\}\right)>0$.
\medskip

Many authors have made contributions in the study of this problem, 
especially for the case $\la=0$. The Palais-Smale condition (PS)
plays a central role when variational methods are applied in the study of problem
\eq{doi}. In this paper, we establish the existence and the multiplicity of 
nontrivial solutions of \eq{doi} with $\la>0$ sufficiently small, 
in a case where the condition (PS) is not assumed even for $\la =0$.
More precisely, we will show 
that there exists at least two weak solutions of \eq{doi}
for $g\not=0$ in an appropriate weighted Sobolev space and $\la>0$ 
small enough. Our proof relies on Ekeland's Variational Principle
\cite{e} and on the Mountain Pass Theorem without the Palais-Smale
condition (in the sense of Brezis and Nirenberg, see \cite{bn}), 
combined with a weighted variant of the
Brezis-Lieb Lemma \cite{bl}.

The natural functional space to study problem
\eq{doi} is $H^1_a(\RR^N)$,  defined as the completion of
$C^{\infty}_0(\RR^N)$ with respect to the norm
\neweq{cinci}
\|u\|=\left(\,\int_{\RR^N}|x|^{-2a}|\nabla u|^2dx\right)
^{1\slash 2}. \endeq It turns out that $H^1_a(\RR^N)$ is a Hilbert
space with respect to the inner product
$$
\langle u,v\rangle=\int_{\RR^N}|x|^{-2a}\nabla u\cdot \nabla
vdx,\quad \forall\:u,v\in H^1_a(\RR^N).
$$
It follows that \eq{patru} holds for all $u\in H^1_a(\RR^N)$.
According to \cite{cw} we have
\neweq{sapte}
\di H^1_a(\RR^N)=\di\overline{C^{\infty}_0\left(\RR^N\setminus
\{0\}\right)}^ {\,\|\cdot\|}, \endeq where $\|\cdot\|$ is 
given by \eq{cinci}.
Let $\|\cdot\|_{-1}$ denote the norm in the dual space $H^{-1}_a(\RR^N)$ 
 of $H^1_a(\RR^N)$. 

Throughout this paper
we suppose that $g\in H^{-1}_a(\RR^N)\setminus \{0\}$.

For an arbitrary open set $\Omega\subset\RR^N$ let $L^p_b(\Omega)$
be the space of all measurable real functions $u$ defined on $\Omega$ 
such that 
$\di\int_{\Omega}|x|^{-bp}\,|u|^pdx$ is finite. By
\eq{patru} it follows that the weighted Sobolev space
$H^1_a(\Omega)$ is continuously
embedded in $L^p_b(\Omega)$.

\medskip
\noindent\textbf{Definition 1.} We say that a function $\,u\in
H^1_a(\RR^N)\,$ is a weak solution of problem \eq{doi} if
$$\di\int_{\RR^N}|x|^{-2a}\nabla u\cdot\nabla v\,dx-\int_{\RR^N}
K(x)|x|^{-bp}\,|u|^{p-2}uv\,dx-\la\int_{\RR^N}g(x)v\,dx=0,$$ for
all $u\in C^{\infty}_0(\RR^N).$
\medskip

 Obviously, the solutions of problem \eq{doi} correspond to
critical points of the energy functional
$$J_{\la}(u)=\frac{1}{2}\int_{\RR^N}|x|^{-2a}|\nabla u|^2\,dx-\frac{1}{p}
\int_{\RR^N}K(x)|x|^{-bp}\,|u|^p\,dx-\la\int_{\RR^N}g(x)u\,dx,$$
where $\,u\in H^1_a(\RR^N)\,$.

\medskip
Our main result is the following.

\begin{teo}\label{t1} Suppose that assumptions $(K1)$, $(K2)$, $(K3)$ are fulfilled
and fix $g\in H^{-1}_a(\RR^N)\setminus\{0\}$. Then there exists
$\,\la_0>0\,$ such that for all $\,\la\in(0,\la_0)$, problem
\eq{doi} has at least two solutions.
\end{teo}

Since the embedding $H^1_a(\RR^N)\hookrightarrow
L^p_b(\RR^N)$ is not compact, the energy functional
$J_\lambda$ fails to satisfy the (PS) condition. Such a failure brings about
difficulty in applying a variational approach to equation \eq{doi}. Furthermore,
since $g\not\equiv0$, then $0$ is no longer a trivial solution of problem \eq{doi} and,
therefore, the Mountain Pass Theorem cannot be applied directly. Using some 
ideas developed in \cite{t}, we obtain the first solution by applying 
Ekeland's Variational Principle. Then, the Mountain Pass Theorem without the 
Palais-Smale condition yields a bounded Palais-Smale sequence whose
weak limit is a critical point of $J_\la$. The proof is concluded by
showing that these two solutions are distinct because
they realize different energy levels.

The paper is organized as follows.
In Section 2 we give some technical results which allow us to give
a variational approach of our main result that we prove in
Section 3. We point out that since the perturbation term $g$ is not assumed to be non-negative
then we can not expect that the distinct solutions given by Theorem~\ref{t1}
are positive. However, if $g\geq 0$ is a non-trivial perturbation, then a straightforward
argument based on the maximum principle implies that the solutions of problem \eq{doi}
are positive.
\medskip

\noindent\textbf{Notations.} Throughout this paper we will denote
by $\,B_R\,$ the open ball in $\,H^1_a(\RR^N)\,$ centered at
origin and having radius $\,R>0.$ We also denote by
$\,\langle\,\cdot\,,\,\cdot\,\rangle\,$ the duality pairing
between $\,H^1_a(\RR^N)\,$ and $\,H^{-1}_a(\RR^N).$ The notations
"$\weak$" and "$\ri$" stand, respectively, for the weak and for the strong
convergence in an arbitrary Banach space.

\section{Auxiliary results}
Define the functionals $\,J_0,I:H^1_a(\RR^N)\rightarrow\RR\,$
by
$$J_0(u)=\frac{1}{2}\int_{\RR^N}|x|^{-2a}|\nabla u|^2\,dx-\frac{1}{p}
\int_{\RR^N}K(x)|x|^{-bp}\,|u|^p\,dx,$$
$$I(u)=\frac{1}{2}\int_{\RR^N}|x|^{-2a}|\nabla u|^2\,dx-\frac{1}{p}
\int_{\RR^N}K_0|x|^{-bp}\,|u|^p\,dx.$$ The
Caffarelli-Kohn-Nirenberg inequality \eq{patru} and the conditions
$(K1)$, $(K2)$ imply that the functionals
$J_{\la}$, $J_0,$ and $I$ are well defined and
$\,J_{\la},J_0,I\in C^1(H^1_a(\RR^N),\RR).$

\medskip
\noindent\textbf{Remark 1.} If
$\,\Omega\subset\RR^N\,$ is a smooth bounded set such that
$\,0\not\in \overline\Omega\,$ then, by the Sobolev
inequality, we have
$$\di\left(\,\int_{\Omega}|x|^{-bp}\,|u|^pdx\right)^{1\slash p}\leq C_1
\left(\,\int_{\Omega}|u|^pdx\right)^{1\slash p}\leq C_2
\left(\,\int_{\Omega}|\nabla u|^2dx\right)^{1\slash 2}\leq C_3
\left(\,\int_{\Omega}|x|^{-2a}|\nabla u|^2dx\right)^{1\slash 2},$$
for all $\,u\in H^1_a(\Omega)\,$. It follows that
$\,H^1_a(\Omega)\,$ is compactly embedded
in $\,L^p_b(\Omega)$.

Remark 1 implies that if $\{u_n\}$ is a sequence that 
converges weakly to some $u_0$ in
$\,H^1_0(\RR^N)$ then $\{u_n\}$ is bounded in
$\,H^1_0(\RR^N)$. Therefore, we can assume (up to a sequence) that
\neweq{opt}
\di u_n\weak u_0 \;\,\mbox{in}\;\,
L^p_{b,\,{\rm loc}}(\RR^N\setminus\{0\}) \quad\mbox{and} \quad
u_n\rightarrow u_0\;\,\mbox{a.e. in}\;\,\RR^N.\\
\endeq

\noindent\textbf{Definition 2.} Let $X$ be a Banach space,
$\,F:X\rightarrow\RR\,$ be a $\,C^1-$functional and $c$ be a real
number. A sequence $\{u_n\}\subset X$ is called a $\,(PS)_c\,$
sequence of $\,F\,$ if $\,F(u_n)\rightarrow c\,$ and
$\,\|F'(u_n)\|_{X^*}\rightarrow 0$.

Our first result shows that if a $\,(PS)_c\,$ sequence of
$\,J_{\la}\,$ is weakly convergent then its limit is a solution of
problem \eq{doi}.

\begin{lemma} \label{l1} Let $\{u_n\}\subset H^1_a(\RR^N)$ be a $\,(PS)_c\,$ 
sequence of
$\,J_{\la}\,$ for some $c\in\RR$. Suppose that $\{u_n\}$ converges
weakly to some $u_0$ in $\,H^1_a(\RR^N)$. Then $\,u_0\,$ is a
solution of problem \eq{doi}.
\end{lemma}
\noindent\proof Let $\,\phi\in C^{\infty}_0(\RR^N\setminus\{0\})\,$ be an
arbitrary function and set $\,\Omega:=\supp\,\phi.$ Since
$\,J'_{\la}(u_n)\rightarrow 0\,$ in $\,H^{-1}_a(\RR^N)\,$ we obtain
$\,\langle J'_{\la}(u_n), \phi\,\rangle\ri 0\,$ as
$\,n\ri\infty\,$, that is,
\neweq{noua}
\di\lim\limits_{n\rightarrow\infty}\left(\,\int_{\Omega}|x|^{-2a}\nabla
u_n \cdot\nabla
\phi\,dx-\int_{\Omega}K(x)|x|^{-bp}\,|u_n|^{p-2}u_n\phi\,dx-
\la\int_{\Omega}g(x)\phi\,dx\right)=0.
\endeq
Since $\,u_n\weak u_0\,$ in $\,H^1_a(\RR^N),$ it follows that
\neweq{zece}
\di\lim\limits_{n\rightarrow\infty}\int_{\Omega}|x|^{-2a}\nabla
u_n \cdot\nabla \phi\,dx=\int_{\Omega}|x|^{-2a}\nabla u_0
\cdot\nabla \phi\,dx.
\endeq

The boundedness of $\{u_n\}$ in $\,H^1_a(\RR^N)\,$ and the
Caffarelli-Kohn-Nirenberg inequality imply that
$\,\{|u_n|^{p-2}u_n\}\,$ is bounded in
$\,L^{p/p-1}_b(\RR^N)$. Since $\,|u_n|^{p-2}u_n\rightarrow
|u_0|^{p-2}u_0\,$ a.e. in $\RR^N$ (which is a consequence of
\eq{opt}), we deduce that $|u_0|^{p-2}u_0$ is the weak
limit in
$L^{p/p-1}_b(\RR^N)$ of the sequence $\{|u_n|^{p-2}u_n\}$. Therefore
\neweq{unspe}
\di\lim\limits_{n\rightarrow\infty}\int_{\Omega}K(x)|x|^{-bp}\,|u_n|^{p-2}u_n
\phi\,dx=\int_{\Omega}K(x)|x|^{-bp}\,|u_0|^{p-2}u_0\phi\,dx.
\endeq
Consequently, relations \eq{noua}, \eq{zece}, and \eq{unspe} yield
$$\di\int_{\Omega}|x|^{-2a}\nabla u_0
\cdot\nabla
\phi\,dx-\int_{\Omega}K(x)|x|^{-bp}\,|u_0|^{p-2}u_0\phi\,dx-
\la\int_{\Omega}g(x)\phi\,dx=0.$$ By virtue of \eq{sapte} we
deduce that the above equality holds for all $\,\phi\in
H^1_a(\RR^N)\,$ which means that $\,J'_{\la}(u_0)=0$. The proof of
our lemma is now complete.
\qed

\bigskip
We now establish a weighted variant of the Brezis-Lieb Lemma (see
\cite{bl}).
\begin{lemma} \label{l2} Let $\{u_n\}$ be a sequence which is weakly convergent to
$\,u_0\,$ in $\,H^1_a(\RR^N)\,$. Then
$$\di\lim\limits_{n\rightarrow\infty}\int_{\RR^N}K(x)|x|^{-bp}\,(|u_n|^p-
|u_n-u_0|^p)\,dx=\int_{\RR^N}K(x)|x|^{-bp}\,|u_0|^p\,dx.$$
\end{lemma}
\noindent\proof Using the boundedness of $\{u_n\}$ in $\,H^1_a(\RR^N)\,$
and the Caffarelli-Kohn-Nirenberg inequality, it follows that the sequence
$\{u_n\}$ is bounded in $\,L^p_b(\RR^N)\,$. Let $\ep>0$
be a positive real number. By $(K1)$ and $(K2)$, we can choose
$\,R_\ep>r_\ep>0\,$ such that
\neweq{doispe}
\di\int_{|x|<r_\ep}K(x)|x|^{-bp}\,|u_0|^p\,dx<\ep,\endeq
and
\neweq{treispe}
\di\int_{|x|>R_\ep}K(x)|x|^{-bp}\,|u_0|^p\,dx<\ep.\endeq
Denote $\,\Omega_\ep=\overline{B(0,R_\ep)}\,\setminus B(0,r_\ep)$. We have
$$\begin{array}{lll}
& \di\left|\;\int_{\RR^N}K(x)|x|^{-bp}\,(|u_n|^p-|u_0|^p-|u_n-u_0|^p)\,dx\right|\\
&
\di\quad\leq\left|\;\int_{\Omega_\ep}K(x)|x|^{-bp}\,(|u_n|^p-|u_0|^p)\,dx\right|
+\int \limits_{\Omega_\ep}K(x)|x|^{-bp}\,|u_n-u_0|^p\,dx\\
& \di\qquad+\int
\limits_{|x|<r_\ep}K(x)|x|^{-bp}\,|u_0|^p\,dx+\left|\;\int_
{|x|<r_\ep}K(x)|x|^{-bp}\,(|u_n|^p-|u_n-u_0|^p)\,dx\right|\\
& \di\qquad+\int
\limits_{|x|>R_\ep}K(x)|x|^{-bp}\,|u_0|^p\,dx+\left|\;\int_
{|x|>R_\ep}K(x)|x|^{-bp}\,(|u_n|^p-|u_n-u_0|^p)\,dx\right|\\
\end{array}$$
By the Lagrange Mean Value Theorem we have
\neweq{paispe}
\di\int_
{|x|<r_\ep}K(x)|x|^{-bp}\,(|u_n|^p-|u_n-u_0|^p)\,dx=p\int_
{|x|<r_\ep}K(x)|x|^{-bp}\,|\theta u_0+(u_n-u_0)|^{p-1}|u_0|\,dx\,,
\endeq where $0<\theta(x)< 1$. Next, we employ the following
elementary inequality: for all $s>0$ there exists a constant $c=c(s)$
such that
$$ (x+y)^s\leq c(x^s+y^s)\quad\mbox{for any}\;\; x,y\in(0,\infty).$$
Then, by H\"{o}lder's inequality and relation \eq{doispe} we deduce that
$$\begin{array}{lll}
& \di \int_{|x|<r_\ep}K(x)|x|^{-bp}\,|\theta u_0+(u_n-u_0)|^{p-1}|u_0|\,dx \\
& \di \qquad\quad\leq
c\!\int_{|x|<r_\ep}K(x)|x|^{-bp}\,(|u_0|^p+|u_n-u_0|^
{p-1}|u_0|)\,dx\\&
\di\qquad\quad=c\int_{|x|<r_\ep}K(x)|x|^{-bp}\,|u_0|^p\,dx
+c\int_{|x|<r_\ep}K(x)|x|^{-bp}\,|u_n-u_0|^{p-1}|u_0|\,dx\\
& \di \qquad\quad \leq c\,\ep+c\left(\,\int_{|x|<r_\ep}
K(x)|x|^{-bp}\,|u_n-u_0|^{p}\,dx\right)^{(p-1)/p}
\left(\,\int_{|x|<r_\ep}K(x)|x|^{-bp}\,|u_0|^p\,dx\right)^{1/p}\\
& \di \qquad\quad \leq c_1\,(\ep+\ep^{1/p}),\\
\end{array}$$
where the constant $\,c_1\,$ is independent of $n$ and $\ep\,$.
Using relation \eq{paispe} we have
\neweq{cincispe}
\di\int
\limits_{|x|<r_\ep}K(x)|x|^{-bp}\,|u_0|^p\,dx+\left|\;\int_{|x|<r_\ep}K(x)|x|^
{-bp}\,(|u_n|^p-|u_n-u_0|^p)\,dx\right|\leq p\,\tilde
c_1\,(\ep+\ep^{1/p}).
\endeq
In a similar manner we obtain
\neweq{saispe}
\di\int
\limits_{|x|>R_\ep}K(x)|x|^{-bp}\,|u_0|^p\,dx+\left|\;\int_{|x|>R_\ep}K(x)|x|^
{-bp}\,(|u_n|^p-|u_n-u_0|^p)\,dx\right|\leq p\,\tilde
c_2\,(\ep+\ep^{1/p}).
\endeq
Since $\,u_n\rightharpoonup u_0$ in $H^1_a(\RR^N),$ relation
\eq{opt} yields
\neweq{saptispe}
\begin{tabular}{ll}
$\di\lim\limits_{n\ri\infty}\int_{\Omega_\ep}K(x)|x|^{-bp}\,(|u_n|^p-|u_0|^p)\,dx=0,$\\
$\di\lim\limits_{n\ri\infty}\int_{\Omega_\ep}K(x)|x|^{-bp}\,|u_n-u_0|^p\,dx=0.$\\
\end{tabular}
\endeq
Now, relations \eq{cincispe}, \eq{saispe}, and \eq{saptispe} yield
$$\di\limsup\limits_{n\ri\infty}\left|\;\int_{\RR^N}K(x)|x|^{-bp}\,(|u_n|^p-|u_0|^p
-|u_n-u_0|^p)\,dx\right|\leq(pC+1)\,(\ep+\ep^{1/p}).$$ Since
$\ep>0$ is arbitrary, it follows that
$$\di\lim\limits_{n\ri\infty}\int_{\RR^N}K(x)|x|^{-bp}\,(|u_n|^p-
|u_n-u_0|^p)\,dx=\int_{\RR^N}K(x)|x|^{-bp}\,|u_0|^p\,dx.$$ This
concludes the proof.
\qed

\begin{lemma} \label{l3} Let $\{v_n\}$ be a sequence which converges weakly to 0 in
$\,H^1_a(\RR^N).$ Then the following properties hold true
$$
\di\lim\limits_{n\ri\infty}[J_{\la}(v_n)-I(v_n)]=0,$$
$$
\di\lim\limits_{n\ri\infty}[\langle
J'_{\la}(v_n),v_n\rangle-\langle I'(v_n),v_n \rangle]=0.$$
\end{lemma}
\noindent\proof A simple computation yields
$$\di J_{\la}(v_n)=I(v_n)-\frac{1}{p}
\int_{\RR^N}(K(x)-K_0)|x|^{-bp}\,|v_n|^p\,dx-\la\int_{\RR^N}g(x)v_n\,dx,$$
$$\di\langle J'_{\la}(v_n),v_n\rangle=\langle I'(v_n),v_n\rangle-
\int_{\RR^N}(K(x)-K_0)|x|^{-bp}\,|v_n|^p\,dx-\la\int_{\RR^N}g(x)v_n\,dx.$$

\noindent Since $\,v_n\weak 0\,$ in $\,H^1_a(\RR^N)$, it follows from the
above equalities that it suffices to prove that
\neweq{douazeci}
\lim_{n\ri\infty}\int_{\RR^N}(K(x)-K_0)|x|^{-bp}\,|v_n|^p\,dx=0.\endeq
Fix $\ep>0.$ By our assumptions $\,(K1)\,$ and
$(K2)$, there exists $\,R_\ep>r_\ep>0\,$ such
that
$$|K(x)-K_0|=K(x)-K_0<\ep \quad\mbox{for \,a.e.} 
\quad x\in \RR^N\setminus\Omega_\ep,$$
where $\Omega_\ep=\overline{B(0,R_\ep)}\setminus B(0,r_\ep)$. Next,
we have
$$\begin{tabular}{ll}
$\quad\di\int_{\RR^N}(K(x)-K_0)|x|^{-bp}\,|v_n|^p\,dx$\\
$\qquad\qquad\di=\int_{\RR^N\setminus\Omega_\ep}
(K(x)-K_0)|x|^{-bp}\,|v_n|^p\,dx+\int_{\Omega_\ep}
(K(x)-K_0)|x|^{-bp}\,|v_n|^p\,dx$\\
$\qquad\qquad\di\leq\ep\int_{\RR^N\setminus\Omega_\ep}|x|^{-bp}\,|v_n|^p\,dx+
(\|K\|_{\infty}-K_0)\int_{\Omega_\ep}|x|^{-bp}\,|v_n|^p\,dx$\\
$\qquad\qquad\di\leq\ep\int_{\RR^N}|x|^{-bp}\,|v_n|^p\,dx+
(\|K\|_{\infty}-K_0)\int_{\Omega_\ep}|x|^{-bp}\,|v_n|^p\,dx.$\\
\end{tabular}$$
Since $v_n\weak 0$ in $H^1_a(\RR^N)$, the Caffarelli-Kohn-Nirenberg
inequality implies
that $\{v_n\}$ is bounded in $L^p_b(\RR^N)$. Moreover, by
 \eq{opt}, it follows that $v_n\ri 0$ in
$L^p_{b,\,{\rm loc}}(\RR^N\setminus\{0\})$. The above relations yield
$$\di\limsup_{n\ri\infty}\int_{\RR^N}(K(x)-K_0)|x|^{-bp}\,|v_n|^p\,dx\leq C\ep$$
for some constant $\,C>0\,$ independent of $n$ and $\ep$. Since
$\,\ep>0\,$ was arbitrarily chosen, we conclude that \eq{douazeci}
holds and the proof of Lemma~\ref{l3} is now complete.
\qed
\medskip

\begin{lemma}\label{l4} There exists $\,\la_1>0\,$ and $\,R=R(\la_1)>0\,$ such
that for all $\la\in(0,\la_1)$, the functional $J_\la$ admits a
$\,(PS)_{c_{0,\la}}\,$ sequence  with
$\,c_{0,\la}=c_{0,\la}(R)=\inf\limits_{u\in\overline{B}_R}J_{\la}(u)$.
Moreover, $\,c_{0,\la}\,$ is achieved by some $\,u_0\in
H^1_a(\RR^N)\,$ with $\,J'_{\la}(u_0)=0\,$.
\end{lemma}
\medskip

\noindent\proof Fix $\la\in(0,1)$. For all $\,u\in
H^1_a(\RR^N)\,$, the assumption $(K_1)$ and the
Caffarelli-Kohn-Nirenberg inequality imply
\medskip

$$\begin{tabular}{ll}
$\di
J_{\la}(u)=\frac{1}{2}\|u\|^2-\frac{1}{p}\int_{\RR^N}K(x)|x|^{-bp}\,|u|^p\,dx-
\la\int_{\RR^N}g(x)u\,dx$\\
$\qquad\;\;\,\di\geq\frac{1}{2}\|u\|^2-\frac{\|K\|_{\infty}}{p}\;C^p_{a,b}\|u\|^p-\la
\|g\|_{-1}\|u\|.$\\
\end{tabular}$$
We now apply the inequality
$\,\di\alpha\beta\leq\frac{\alpha^2+\beta^2}{2},$ for
any $\alpha,\beta\geq 0.$ Hence
\neweq{douazecisiunu}
\di
J_{\la}(u)\geq\frac{1-\la}{2}\|u\|^2-\frac{\|K\|_{\infty}}{p}\:C^p_{a,b}\,\|u\|^p-
\frac{\la}{2}\|g\|^2_{-1}.\endeq Since $\,p>2\,$ and the right
side of \eq{douazecisiunu} is a decreasing function on $\,\la\,$,
we find $\,\la_1>0\,$ and
$R=R(\la_1)>0,\,\delta=\delta(\la_1)>0\,$ such that
\neweq{douazecibis}
J_{\la}(u)\geq -\frac{\la}{2}\|g\|^2_{-1},\quad\mbox{ for all }
\;u\in\overline B_R\;\; \mbox{ and } \la\in(0,\la_1)
\end{equation}
and
\neweq{douazecibisbis}
J_{\la}(u)\geq\delta>0,\quad\mbox{ for all } \;u\in\partial
B_R\;\; \mbox{ and } \la\in(0,\la_1).
\end{equation}
For instance, we can take
$$\di \la_1:=\min\left\{\frac{1}{2},\frac{1}{2\|g\|^2_{-1}}
\left(\frac{1}{2}-\frac{1}{p}\right)r_0^2\right\},$$
$$\di r_0:=\left[\frac{1}{2\|K\|_{\infty}C^p_{a,b}}\right]^{1/(p-2)},\quad
 R:=\left[\frac{1-\la_1}{\|K\|_{\infty}C^p_{a,b}}\right]^{1/(p-2)}$$
and
$$\quad\delta(\la_1):=\frac{\la_1}{2}\|g\|^2_{-1}.$$
Using now the estimate \eq{douazecisiunu}, we easily deduce
\eq{douazecibis} and
\eq{douazecibisbis}.

Next, we define
$\,c_{0,\la}:=c_{0,\la}(R)=\inf\{J_\la(u)\,;\,u\in\overline{B}_R
\}$. We first note that $\,c_{0,\la}\leq J_\la(0)=0\,$. The set
$\,\overline{B}_R\,$ becomes a complet metric space with respect
to the distance
$$\mbox{dist}(u,v)=\|u-v\|\,,\qquad\mbox{for any}\;\;u,v\in\overline{B}_R.$$
The functional $\,J_\la\,$ is lower semi-continuous and bounded
from below on $\,\overline{B}_R\,$. Then, by Ekeland's Variational
Principle \cite[Theorem 1.1]{e}, for any positive integer $\,n\,$
there exists $u_n$ such that
\neweq{douazecisidoi}
c_{0,\la}\leq J_\la(u_n)\leq c_{0,\la}+\frac{1}{n}\endeq and
\neweq{douazecisitrei}
J_\la(w)\geq J_\la(u_n)-\frac{1}{n}\|u_n-w\|\quad\mbox{for
all}\;\;w\in\overline{B}_R.
\endeq
We first show that $\,\|u_n\|<R\,$ for $n$ large enough. Indeed,
if not, then $\,\|u_n\|=R\,$ for infinitely many $\,n,$ and so (up
to a subsequence) we can assume  that $\,\|u_n\|=R\,$ for all
$\,n\geq 1$. It follows that $\,J_\la(u_n)\geq\delta>0\,$. Using
\eq{douazecisidoi} and letting $\,n\ri\infty\,$, we have $\,0\geq
c_{0,\la}\geq\delta>0$, which is a contradiction.

We now claim that $\,J'_\la(u_n)\ri 0\,$ in $\,H^{-1}_a(\RR^N)\,$.
Fix $\,u\in H^1_a(\RR^N)\,$ with $\,\|u\|=1\,$ and let
$\,w_n=u_n+tu$. For some fixed $\,n\,$, we have $\|w_n\|\leq
\|u_n\|+t<R\,$ if $\,t>0\,$ is small enough. Then relation

\eq{douazecisitrei} yields
$$J_\la(u_n+tu)\geq J_\la(u_n)-\frac{t}{n}\|u\|\,,$$
that is,
$$\frac{J_\la(u_n+tu)-J_\la(u_n)}{t}\geq -\frac{1}{n}\|u\|=-\frac{1}{n}.$$
Letting $\,t\searrow 0\,$ it follows that $\di\,\langle
J'_\la(u_n),u\rangle\geq- \frac{1}{n}$. Arguing in a similar way
for $\,t\nearrow 0,$ we obtain $\di\,\langle J'_\la(u_n),
u\rangle\leq\frac{1}{n}$. Since $\,u\in H^1_a(\RR^N)\,$ with
$\,\|u\|=1\,$ has been arbitrarily chosen, we have
$$\|J'_\la(u_n)\|=\sup\limits_{u\in H^1_a(\RR^N),\\\|u\|=1}|\langle J'_\la(u_n),u\rangle|
\leq\frac{1}{n}\ri 0\quad\mbox{as}\;\;n\ri\infty.$$
\medskip

We  have proved the existence of a $\,(PS)_{c_{0,\la}}\,$
sequence, i.e., a sequence $\,\{u_n\}\subset H^1_a(\RR^N)\,$ with
\neweq{douazecisipatru}
J_\la(u_n)\ri c_{0,\la}\quad\mbox{and}\quad J'_\la(u_n)\ri
0\;\;\mbox{in}\;\;H^1_a(\RR^N).
\endeq
Since $\,\|u_n\|\leq R\,$, it follows that $\{u_n\}$ converges weakly
(up to a subsequence)
 in $\,H^1_a(\RR^N)\,$ to some $\,u_0$. Moreover, 
 relations \eq{opt}, \eq{douazecisipatru}, and Remark
1 yield
\neweq{douazecisicinci}
u_n\weak u_0 \quad\mbox{in}\;\;H^1_a(\RR^N),\qquad u_n\ri u_0
\quad\mbox{a.e. in}\;\; \RR^N\endeq and
\neweq{douazecisisase}
J'_\la(u_0)=0.\endeq

Next, we prove that $J_\la(u_0)=c_{0,\la}$. Using relations
\eq{douazecisipatru} and \eq{douazecisicinci} we have
$$\di o(1)=\langle J'_\la(u_n),u_n\rangle=\int_{\RR^N}|x|^{-2a}|\nabla u_n|^2\,dx-
\int_{\RR^N}K(x)|x|^{-bp}\,|u_n|^p\,dx-\la\int_{\RR^N}g(x)u_n\,dx.$$
Therefore
$$\di J_\la(u_n)=\left(\frac{1}{2}-\frac{1}{p}\right)\int_{\RR^N}K(x)|x|^{-bp}|u_n|
^p\,dx-\frac{\la}{2}\int_{\RR^N}g(x)u_n\,dx+o(1).$$ 
Hence
$$\di J_\la(u_0)=\left(\frac{1}{2}-\frac{1}{p}\right)\int_{\RR^N}K(x)|x|^{-bp}|u_0|
^p\,dx-\frac{\la}{2}\int_{\RR^N}g(x)u_0\,dx+o(1).$$ 
Fatou's Lemma
and relations \eq{douazecisipatru}, \eq{douazecisicinci},
\eq{douazecisisase} imply
$$ c_{0,\la}=\liminf\limits_{n\ri\infty}\di J_\la(u_n)
\geq\left(\frac{1}{2}-\frac{1}{p}\right)
\int_{\RR^N}K(x)|x|^{-bp}|u_0|^p\,dx-\frac{\la}{2}\int_{\RR^N}g(x)u_0\,dx
\di=J_\la(u_0).
$$
Thus, $c_{0,\la}\geq J_\la(u_0).$ On the other
hand, since $u_0\in \overline B_R,$ we deduce that $J_\la(u_0)\geq
c_{0,\la},$ so $J_\la(u_0)=c_{0,\la}.$ This concludes the proof of Lemma \ref{l4}.
\qed

\section{Proof of Theorem \ref{t1}}

Define
$${\cal S}=\{u\in H^1_a(\RR^N)\setminus\{0\}\,;\,\langle I'(u),u\rangle=0\}.$$
We claim that ${\cal S}\not=\emptyset$. For this purpose we fix
$u\in H^1_a(\RR^N)\setminus\{0\}$ and set, for any $\la>0,$
$$\di\Psi(\la)=\langle I'(\la u),\la u\rangle=\la^2\int_{\RR^N}|x|^{-2a}|\nabla u|
^2\,dx-\la^p\int_{\RR^N}K_0|x|^{-bp}\,|u|^p\,dx.$$ Since $p>2$, it
follows that $\Psi(\la)<0$ for $\la$ large enough and
$\Psi(\la)>0$ for $\la$ sufficiently close to the origin. So, there
exists $\la>0$ such that $\Psi(\la)=0,$ that is, $\la u\in {\cal
S}.$
\medskip

\begin{prop}\label{p1}
Let $I_{\infty}:=\inf\{\,I(u)\:;\,u\in {\cal S}\}.$ Then there
exists $\bar u\in H^1_a(\RR^N)$ such that
\neweq{douazecisisapte}
I_{\infty}=I(\bar u)=\sup\limits_{t\geq 0}I(t\bar u).
\endeq
\end{prop}

\noindent\proof For some fixed $\phi\in H^1_a(\RR^N)\setminus \{0\}$ denote
$$\di f(t)=I(t\phi)=\frac{t^2}{2}\int_{\RR^N}|x|^{-2a}|\nabla \phi|^2\,dx-
\frac{K_0}{p}\,t^p\int_{\RR^N}|x|^{-bp}\,|\phi|^p\,dx.$$ We have
$$\di f'(t)=t\int_{\RR^N}|x|^{-2a}|\nabla \phi|^2\,dx-
K_0\,t^{p-1}\int_{\RR^N}|x|^{-bp}\,|\phi|^p\,dx.$$ 
Then $f$ attains its maximum at
$$\di t_0=t_0(\phi):=\left\{ \frac{\di\int_{\RR^N}|x|^{-2a}|\nabla \phi|^2\,dx}
{\di \int_{\RR^N}K_0|x|^{-bp}|\phi|^p\,dx} \right\}^{1\slash
(p-2)}.$$ Hence
$$\di f(t_0)=I(t_0\phi)=\sup\limits_{t\geq 0}I(t\phi)=\left(\,\frac{1}{2}-
\frac{1}{p}\right)\left\{ \frac{\di\int_{\RR^N}|x|^{-2a}|\nabla
\phi|^2\,dx}
{\di\left(\int_{\RR^N}K_0|x|^{-bp}|\phi|^p\,dx\right)^{2/p}}
\right\}^{p\slash (p-2)}.$$ It follows that
\neweq{douazecisiopt}
\di\inf\limits_{\phi\in
H^1_0(\RR^N)\setminus\{0\}}\,\sup\limits_{t\geq 0}I(t\phi)=
\left(\,\frac{1}{2}-\frac{1}{p}\right)\left[S(a,b)\right]^{p/(p-2)},\endeq
where
\neweq{douazecisinoua}
\di S(a,b)=\inf\limits_{\phi\in H^1_0(\RR^N)\setminus\{0\}}
\left\{ \frac{\di\int_{\RR^N}|x|^{-2a}|\nabla \phi|^2\,dx}
{\di\left(\int_{\RR^N}K_0|x|^{-bp}|\phi|^p\,dx\right)^{2/p}}
\right\}.\endeq We now easily observe that for every $u\in {\cal
S}$ we have $\,t_0(u)=1$, so, by \eq{douazecisiopt} it follows
that
\neweq{treizeci}
I(u)=\sup\limits_{t\geq 0}I(tu) \quad\mbox{for all}\;\;u\in {\cal
S}. \endeq

According to \cite[Theorems 1.2, 7.2, 7.6]{cw}, the infimum
in \eq{douazecisinoua} is achieved by a function $U\in
H^1_a(\RR^N)$ such that $\int_{\RR^N}K_0|x|^{-bp}|U|^pdx=1.$
Letting $\,\bar u=\left[S(a,b)\right]^{1/(p-2)}U$, we see that
$\bar u\in {\cal S}$ and
\neweq{treizecisiunu}
I(\bar
u)=\left(\frac{1}{2}-\frac{1}{p}\right)\left[S(a,b)\right]^{p/(p-2)}.
\endeq
Relations \eq{treizeci} and \eq{treizecisiunu} yield
$$I_{\infty}=\inf\limits_{u\in {\cal S}}I(u)=
\inf\limits_{u\in {\cal S}}\sup\limits_{t\geq 0}I(tu)\geq
\inf\limits_{u\in H^1_0(\RR^N)\setminus\{0\}}\,\sup\limits_{t\geq
0}I(tu)=
\left(\frac{1}{2}-\frac{1}{p}\right)\left[S(a,b)\right]^{p/(p-2)}=I(\bar
u),$$ which concludes our proof.
\qed

\medskip
\begin{prop}\label{p2}
Assume that $\{u_n\}$ is a $(PS)_c$ sequence of $J_\la$ which is
weakly convergent in $H^1_a(\RR^N)$ to some $u_0$. Then the
following alternative holds: either
$\{u_n\}$ converges strongly in $H^1_a(\RR^N),$ or $c\geq
J_\la(u_0)+I_{\infty}.$
\end{prop}

\noindent\proof Since $\{u_n\}$ is a $(PS)_c$ sequence and $u_n\weak u_0$
in $\,H^1_a(\RR^N)$ we have
\neweq{treizecisidoi}
J_\la(u_n)=c+o(1) \quad\mbox{and}\quad\langle
J'_\la(u_n),u_n\rangle=o(1).\endeq Denote $v_n=u_n-u_0$. It
follows that $\,v_n\weak 0\,$ in $\,H^1_a(\RR^N)$ which implies
$$\di\lim\limits_{n\ri\infty}\int_{\RR^N}|x|^{-2a}\nabla v_n\cdot\nabla u_0\,dx=0,$$
$$\di\lim\limits_{n\ri\infty}\int_{\RR^N}g(x)v_n\,dx=0.$$
The above relations imply
\neweq{treizecisitrei}
\begin{tabular}{ll}
$\|u_n\|^2=\|u_0\|^2+\|v_n\|^2+o(1)$\\
$J_\la(v_n)=J_0(v_n)+o(1)$\\
\end{tabular}\endeq
Using Lemmas \ref{l1}-\ref{l3} and relations \eq{treizecisidoi},
\eq{treizecisitrei} we deduce that
\neweq{treizecisipatru}
o(1)+c=J_\la(u_n)=J_\la(u_0)+J_\la(v_n)+o(1)=J_\la(u_0)+I(v_n)+o(1),\endeq
\neweq{treizecisicinci}
o(1)=\langle J_\la'(u_n),u_n\rangle=\langle
J'_\la(u_0),u_0\rangle+\langle J'_\la(v_n),v_n\rangle+o(1)=\langle
I'(v_n),v_n\rangle+o(1).\endeq 
If $v_n\ri 0$ in $H^1_a(\RR^N)\,$
then
$u_n\ri u_0$ in $H^1_a(\RR^N)$. It follows that 
$ J_\la(u_0)=\lim\limits_{n\ri\infty}J_\la(u_n).$
If $v_n\not\ri 0$ in $H^1_a(\RR^N),$ using the fact that $v_n\weak
0$ in $H^1_a(\RR^N),$ we can asume that $\|v_n\|\ri l>0.$

 By virtue of \eq{treizecisipatru}, it remains only to
show that $I(v_n) \geq I_{\infty}+o(1).$ Taking $t>0$ we have
$$\di\langle I'(tv_n),tv_n\rangle=t^2\int_{\RR^N}|x|^{-2a}|\nabla v_n|^2\,dx-
t^p\,K_0\int_{\RR^N}|x|^{-bp}\,|v_n|^p\,dx.$$ If we prove the
existence of a sequence $\{t_n\}\subset(0,\infty)$ with $t_n\ri 1$
and $\langle I'(t_nv_n),t_nv_n\rangle=0,$ then $t_nv_n\in {\cal
S}.$ This implies that
$$\begin{tabular}{ll}
$\di I(v_n)=I(t_nv_n)+\frac{1-t_n^2}{2}\|v_n\|^2-
\frac{1-t_n^p}{p}K_0\int_{\RR^N}|x|^{-bp}\,|v_n|^p\,dx$\\
$\qquad\;\,\,=I(t_nv_n)+o(1)\geq I_{\infty}+o(1),$\\
\end{tabular}$$
and the conclusion follows. For this purpose, we denote
$$\di\alpha_n=\int_{\RR^N}|x|^{-2a}|\nabla v_n|^2\,dx=\|v_n\|^2\geq 0,$$
$$\di\beta_n=K_0\int_{\RR^N}|x|^{-bp}|v_n|^p\,dx\geq 0,$$
$$\mu_n=\alpha_n-\beta_n.$$
From \eq{treizecisicinci} it follows that $\mu_n=\langle
I'(v_n),v_n\rangle\ri 0$ as $n\ri\infty.$ If $\mu_n=0,$ then we
take $t_n=1.$ We next assume that $\mu_n\neq 0.$ Let
$\delta\in\RR$ with $|\delta|>0$ sufficiently small and
$t=1+\delta.$ Then
$$\begin{tabular}{ll}
$\langle
I'(tv_n),tv_n\rangle$&$\di=(1+\delta)^2\alpha_n-(1+\delta)^p\beta_n=(1+\delta)^2
\alpha_n-(1+\delta)^p(\alpha_n-\mu_n)$\\
&$\di=\alpha_n(2\delta-p\delta+o(\delta))+(1+\delta)^p\mu_n$\\
&$\di=\alpha_n(2-p)\delta+\alpha_no(\delta)+(1+\delta)^p\mu_n.$\\
\end{tabular}$$
Since $p>2,$ $\alpha_n\ri l^2>0$ and $\mu_n\ri 0,$ for $n$ large
enough we can define $\di
\delta^+_n=\frac{2|\mu_n|}{\alpha_n(p-2)}$ and $\di\delta^-_n=
-\frac{2|\mu_n|}{\alpha_n(p-2)}.$ It follows that
$$\di \delta^+_n\searrow 0\quad\mbox{and}\quad \langle I'((1+\delta^+_n)v_n),
(1+\delta^+_n)v_n\rangle<0,$$
$$\di \delta^-_n\nearrow 0\quad\mbox{and}\quad \langle I'((1+\delta^-_n)v_n),
(1+\delta^-_n)v_n\rangle<0.$$ From the above relations we deduce
the existence of some
 $t_n\in(1+\delta^-_n,1+\delta^+_n)$ such that
$t_n\ri 1$ and $\langle I'(t_nv_n),t_nv_n\rangle=0.$ This
concludes the proof.
\qed

\medskip
We now fix $\bar u\in H^1_a(\RR^N)$ such that \eq{douazecisisapte}
holds. Since $p>2,$ there exists $\bar t$ such that
$$\begin{tabular}{ll}
$\;\,I(t\bar u)<0\quad\mbox{for all}\quad t>\bar t,$\\
$J_\la(t\bar u)<0\quad\mbox{for all}\quad t>\bar t\;\;\mbox{and}\;\;\la>0.$\\
\end{tabular}$$
Set
\neweq{treizecisisase}
\di{\cal P}=\{\gamma\in C([0,1],H^1_a(\RR^N))\,;\,\gamma(0)=0,\,
\gamma(1)=\bar t\bar u\},
\endeq
\neweq{treizecisisapte}
\di c_g=\inf\limits_{\gamma\in{\cal
P}}\sup\limits_{u\in\gamma}\,J_\la(u).\endeq
\medskip

\begin{prop}\label{l333}
There exists $\la_0>0,$ $R_0=R_0(\la_0)>0,$
$\delta_0=\delta_0(\la_0)>0$ such that $J_\la\mid_{\partial
B_{R_0}}\geq \delta_0$ and $c_g<c_{0,\la}+I_{\infty}$ for all
$\la\in(0,\la_0),$ where
$c_{0,\la}=\inf\limits_{u\in\overline{B}_{R_0}}J_\la(u).$
\end{prop}

\noindent\proof By our hypothesis $(K3)$ and the definition of $I$ we can
assume that
$$J_0(t\bar u)<I(t\bar u)\quad\mbox{for all}\;\;t>0.$$
An elementary computation implies the existence of some $t_0\in(0,\bar t)$
such that
$$\sup\limits_{t\geq 0}J_0(t\bar u)=J_0(t_0\bar u)<I(t_0\bar u)\leq
\sup\limits_{t\geq 0}I(t\bar u)=I_{\infty}.$$ 
So,we can choose
$\ep_0\in(0,1)$ such that
\neweq{treizecisiopt}
\sup\limits_{t\geq 0}J_0(t\bar u)<I_{\infty}-\ep_0.
\endeq
Set
\neweq {treizecisinoua}
\di\la_0:=\min\left\{\la_1,\,\frac{\ep_0}{2\bar t\,\|\bar
u\|\,\|g\|_{-1}},\, \frac{\ep_0}{2\|g\|^2_{-1}} \right\}.\endeq
Applying Lemma 4 it follows that there exists $R_0=R_0(\la_0)>0$
such that for all $\la\in(0,\la_0)$ the conclusion of Lemma 4
holds. Moreover, by virtue of its proof, there exists
$\delta_0=\delta(\la_0)>0$ such that
$J_\la\mid_{\partial\overline{B}_{R_0}}\geq\delta_0.$ Then relations
\eq{treizecisinoua} and \eq{douazecibis} yield
\neweq{patruzeci}
c_{0,\la}=\inf\limits_{u\in\overline{B}_{R_0}}J_\la(u)\geq-\frac{\la}{2}\|g\|^2_{-1}>
-\frac{\ep_0}{2},\quad\mbox{for all}\;\;\la\in(0,\la_0).\endeq For
$u\in\gamma_0=\{t\bar t\bar u\,;\,0\leq t\leq 1\}\in{\cal P}$ we
have
$$\di |J_\la(u)-J_0(u)|=\la\left|\int_{\RR^N}g(x)u\,dx\right|\leq\la\,\bar t\,\|\bar u\|\,
\|g\|_{-1}\leq \frac{\ep_0}{2}\quad\mbox{for
all}\;\;\la\in(0,\la_0).$$ 
Therefore
\neweq{patruzecisiunu}
J_\la(u)\leq J_0(u)+\frac{\ep_0}{2}, \quad\mbox{for
all}\;\;\la\in(0,\la_0).\endeq 
Using relations \eq{treizecisiopt},
\eq{patruzeci} and \eq{patruzecisiunu} we obtain
$$\begin{tabular}{ll}
$\di c_g=\inf\limits_{\gamma\in{\cal
P}}\sup\limits_{u\in\gamma}\,J_\la(u)\leq
\sup\limits_{u\in\gamma_0}\,J_\la(u)$\\
$\di\quad\,\leq\sup\limits_{u\in\gamma_0}\,J_0(u)+\frac{\ep_0}{2}\leq\sup\limits_{t\geq
0}
J_0(t\bar u)+\frac{\ep_0}{2}<I_{\infty}-\frac{\ep_0}{2}<I_{\infty}+c_{0,\la}.$\\
\end{tabular}$$
This completes the proof.
\qed

\bigskip
\textbf{Proof of Theorem \ref{t1} concluded.} Consider $\,R_0>0$,
$\,\delta_0>0\,$ given by Proposition \ref{l333}. In view of
its proof, we deduce that for all $\,\la\in(0,\la_0)\,$ the
conclusion of Lemma 4 holds. Therefore, we obtain the existence of
a solution $\,u_0\,$ of problem \eq{doi} such that
$\,J_\la(u_0)=c_{0,\la}$.

On the other hand, applying the Mountain Pass Theorem without the
Palais-Smale condition (see \cite[ Theorem 2.2]{bn}), it follows that
there exists a $(PS)_{c_g}$ sequence $\{u_n\}$ of $J_\la,$ that is,
$$\di J_\la(u_n)=c_g+o(1)\quad\mbox{and}\quad J'_\la(u_n)\ri 0\;\;\mbox{in}\;\;H^{-1}_a
(\RR^N).$$ Therefore
$$\begin{tabular}{ll}
$\di c_g+o(1)+\frac{1}{p}\|J'_\la(u_n)\|_{-1}\|u_n\|\geq
J_\la(u_n)-\frac{1}{p}\,\langle
J'_\la(u_n),u_n\rangle$\\
$\di\hspace{5,05cm}\geq\left(\frac{1}{2}-\frac{1}{p}\right)\|u_n\|^2-\la\left(1-\frac{1}
{p}\right)\|g\|_{-1}\|u_n\|.$\\
\end{tabular}$$
This above inequality shows that $\{u_n\}$ is  bounded
in $H^1_a(\RR^N).$ Thus we can assume (up to a subsequence) that
$u_n\weak u_1$ in $H^1_a(\RR^N).$ By Lemma 1 it follows that $u_1$
is a weak solution of problem \eq{doi}.

We claim that $u_0\neq u_1.$ Indeed, by Proposition 2, the 
following alternative holds: either
$u_n\ri u_1$ in $H^1_a(\RR^N)$, which gives
$$\di J_\la(u_1)=\lim\limits_{n\ri\infty}J_\la(u_n)=c_g>0\geq c_{0,\la}=J_\la(u_0)$$
and the conclusion follows; or
$$c_g=\lim\limits_{n\ri\infty}J_\la(u_n)\geq J_\la(u_1)+I_{\infty}.$$
In the last case, if we suppose that $u_1=u_0$ then
$J_\la(u_1)=J_\la(u_0)=c_{0,\la}$ and so $c_g\geq
c_{0,\la}+I_{\infty},$ which contradicts Proposition \ref{l333}.
The proof of Theorem 1 is now complete.
\qed

\end{document}